\documentclass[12pt,twoside]{article}

\usepackage{geometry} 
\usepackage{graphicx} 

\usepackage{booktabs} 
\usepackage{array} 
\usepackage{paralist} 
\usepackage{verbatim} 
\usepackage{subfig} 

\usepackage{amssymb}
\usepackage{amsmath}
\usepackage{fancyref}
\usepackage{amsthm}
\usepackage{cases}
\usepackage{hyperref} 

\usepackage{fancyhdr} 
\usepackage{sectsty}
\allsectionsfont{\sffamily\mdseries\upshape} 

\usepackage[nottoc,notlof,notlot]{tocbibind} 
\usepackage[titles,subfigure]{tocloft} 

\numberwithin{equation}{section}

\title{\Large{\bf New Inequalities and Applications}}
\author{\bf {Daiyuan Zhang}}
\begin{document}
\maketitle
\centerline{College of Computer}
\centerline{Nanjing University of Posts and Telecommunications}
\centerline{Nanjing, P.R. China}
\centerline{dyzhang@njupt.edu.cn, zhangdaiyuan2011@sina.com}

\newtheorem{Theorem}{\quad Theorem}[section]
\newtheorem{Definition}[Theorem]{\quad Definition}
\newtheorem{Corollary}[Theorem]{\quad Corollary}
\newtheorem{Lemma}[Theorem]{\quad Lemma}
\newtheorem{Example}[Theorem]{\quad Example}

\centerline{}

\begin{abstract}
This paper presents some new inequalities, the most important of which is the inequality given in Theorem 2.1. It can solve a class of inequalities by a unified method. An important application of the inequality given in Theorem 2.1 is to derive another new general form of inequality. The famous Nesbitt's inequality is a special case of this general form of inequality when $n=3$.

The new inequality in Theorem 2.1 proposed in this paper is easy to use and expand, and many new inequalities can be derived and obtained by direct calculation, so it has a wide range of applications. 

Many known inequalities can also be directly calculated by the inequalities proposed in this paper, and the calculation is simple and convenient. 
\end{abstract}



{\footnotesize Copyright $\copyright$ This is my original work. Prohibit any copying or plagiarism.}

\section{Introduction}\label{sec1}
As we all know, inequalities are widely used in mathematics and natural science. Compared with other branches of mathematics, inequality develops relatively late. It is believed that Hardy, Littlewood and Polya are the first ones to systematically discuss inequalities (see Ref. [1]). In the following decades, documents [2] and [3] also systematically discuss inequalities.

However, the proof and application of inequalities often require relatively high skills. To exaggerate, the problem of inequalities is basically a problem-by-problem solution, which is too "concrete analysis of specific problems", and lacks a unified way of thinking and systematic methods. Because of the relative independence of inequality methods and the seemingly lack of intrinsic correlation, researchers and applicants will often encounter a lot of confusion.

Obviously, it is impossible to prove and apply all inequalities in one way, but people hope to provide as systematic a method as possible, so that this method can solve a class of problems or a large class of problems, so as to avoid the unrelated "concrete analysis of specific problems" as far as possible.

This paper makes an attempt to link some seemingly unrelated inequalities and deal with them in a unified way. For example, the following inequalities seem to be independent of each other (where ${{x}_{i}}\in{{\mathbf{R}}^{+}}$, if and only if each ${{x}_{i}}$ is equal to each other, the following equalities hold):
\[\frac{{{x}_{1}}}{{{x}_{2}}+{{x}_{3}}}+\frac{{{x}_{2}}}{{{x}_{1}}+{{x}_{3}}}+\frac{{{x}_{3}}}{{{x}_{1}}+{{x}_{2}}}\ge \frac{3}{2};\]
 \[\frac{1}{{{x}_{2}}+{{x}_{3}}}+\frac{1}{{{x}_{3}}+{{x}_{1}}}+\frac{1}{{{x}_{1}}+{{x}_{2}}}\ge \frac{9}{2\left( {{x}_{1}}+{{x}_{2}}+{{x}_{3}} \right)};\]
\[\frac{1}{1-x_{1}^{2}}+\frac{1}{1-x_{2}^{2}}\ge \frac{4}{1+2{{x}_{1}}{{x}_{2}}}, {{x}_{1}}+{{x}_{2}}=1;\]
\[\begin{aligned}
& \frac{1}{{{x}_{1}}+{{x}_{2}}}+\frac{1}{{{x}_{1}}+{{x}_{3}}}+\frac{1}{{{x}_{1}}+{{x}_{4}}}+\frac{1}{{{x}_{2}}+{{x}_{3}}}+\frac{1}{{{x}_{2}}+{{x}_{4}}}+\frac{1}{{{x}_{3}}+{{x}_{4}}} \\ 
& \ge \frac{12}{{{x}_{1}}+{{x}_{2}}+{{x}_{3}}+{{x}_{4}}} \\ 
\end{aligned}\]. 

However, through this study, we know that they all originate from one inequality (or that these inequalities have the same "father"), which is the inequality I put forward in this paper.

Among those inequalities mentioned above, the first two are well known, and the first one is the famous Nesbitt's inequality. However, the last two inequalities may be rarely known.

\section{New Inequalities}\label{sec2}
I have proved the following theorem. 

\begin{Theorem}\label{thm2.1}
Let $S$ be a constant, $S-{{x}_{i}}>0$, $i=1,2,\cdots ,n$, $\bar{x}={\sum\limits_{i=1}^{n}{{{x}_{i}}}}/{n}\;$, then the following inequality holds
	\begin{equation}\label{eq2.1}
	\sum\limits_{i=1}^{n}{\frac{1}{S-{{x}_{i}}}}\ge \frac{n}{S-\bar{x}}
	\end{equation} 
Equality holds iff ${{x}_{1}}={{x}_{2}}=\cdots ={{x}_{n}}$. 
\end{Theorem}

\subsection{Generalization of Nesbitt's inequality}\label{sec2.1}
Inequality (2.1) given in theorem 2.1 is very important and has many applications. In this section, I first generalize Nesbitt's inequality to a general form by using inequality (2.1) that I proposed. Other applications are given later in section 3.

It is well known that the following inequality holds for any given positive real number ${{x}_{1}}$, ${{x}_{2}}$ and ${{x}_{3}}$ 
\[\frac{{{x}_{1}}}{{{x}_{2}}+{{x}_{3}}}+\frac{{{x}_{2}}}{{{x}_{1}}+{{x}_{3}}}+\frac{{{x}_{3}}}{{{x}_{1}}+{{x}_{2}}}\ge \frac{3}{2}\]
Equality holds iff ${{x}_{1}}={{x}_{2}}={{x}_{3}}$.

The above inequality is called Nesbitt's inequality. Nesbitt's inequality has also been discussed in some literatures (e.g. [4], [5], [6], [7], [8], etc.), but they are basically limited to Nesbitt's inequality of three quantities, four quantities and six quantities. Although the proof of Nesbitt's inequality is given in those literatures (e.g. [4], [5], [6], [7], [8], etc.), the proof methods are individualized, lack generality, and are not easy to generalize to the general form.

If any given $n$ positive real numbers  ${{x}_{1}}$,${{x}_{2}}$,…,${{x}_{n}}$, How to determine the lower bound of the following expression?
\[\frac{{{x}_{1}}}{{{x}_{2}}+{{x}_{3}}+\cdots +{{x}_{n}}}+\frac{{{x}_{2}}}{{{x}_{1}}+{{x}_{3}}+\cdots +{{x}_{n}}}+\cdots +\frac{{{x}_{n}}}{{{x}_{1}}+{{x}_{2}}+\cdots +{{x}_{n-1}}}\] 

This requires that the Nesbitt's inequality be extended to the general form of $n$ variables. 
The result is given in this paper. Moreover, under certain conditions, the upper bound of the above expression is also given, see theorem 2.3. 

\begin{Theorem}\label{thm2.2}
	Let ${{x}_{i}}\in{{\mathbf{R}}^{+}}$, $i=1,2,\cdots, n$, then the following inequality holds
	\begin{equation}\label{eq2.2}
	\sum\limits_{i=1}^{n}{\frac{{{x}_{i}}}{\sum\limits_{j=1}^{n}{{{x}_{j}}}-{{x}_{i}}}}\ge \frac{n}{n-1}
	\end{equation} 
	Equality holds iff ${{x}_{1}}={{x}_{2}}=\cdots ={{x}_{n}}$. 
\end{Theorem}

{\bf Proof.}
Let 
\begin{equation}\label{eq2.3}
	S=\sum\limits_{j=1}^{n}{{{x}_{j}}}
\end{equation}

According to inequality (2.1), we have
\[\sum\limits_{i=1}^{n}{\frac{\sum\limits_{j=1}^{n}{{{x}_{j}}}}{\sum\limits_{j=1}^{n}{{{x}_{j}}}-{{x}_{i}}}}\ge \frac{{{n}^{2}}}{n-1}\]	
\[\sum\limits_{i=1}^{n}{\frac{\sum\limits_{j=1}^{n}{{{x}_{j}}}-{{x}_{i}}+{{x}_{i}}}{\sum\limits_{j=1}^{n}{{{x}_{j}}}-{{x}_{i}}}}\ge \frac{{{n}^{2}}}{n-1}\]	
\[\sum\limits_{i=1}^{n}{\left( 1+\frac{{{x}_{i}}}{\sum\limits_{j=1}^{n}{{{x}_{j}}}-{{x}_{i}}} \right)}\ge \frac{{{n}^{2}}}{n-1}\]	
i.e.
\[\sum\limits_{i=1}^{n}{1}+\sum\limits_{i=1}^{n}{\frac{{{x}_{i}}}{\sum\limits_{j=1}^{n}{{{x}_{j}}}-{{x}_{i}}}}\ge \frac{{{n}^{2}}}{n-1}\] 
\[n+\sum\limits_{i=1}^{n}{\frac{{{x}_{i}}}{\sum\limits_{j=1}^{n}{{{x}_{j}}}-{{x}_{i}}}}\ge \frac{{{n}^{2}}}{n-1}\]
i.e.
\begin{equation}\label{eq.2.4}
	\sum\limits_{i=1}^{n}{\frac{{{x}_{i}}}{\sum\limits_{j=1}^{n}{{{x}_{j}}}-{{x}_{i}}}}\ge \frac{n}{n-1}
\end{equation}
\qed

Let $n=3$, From inequality (2.4), we get Nesbitt's inequality. Nesbitt's inequality is a special case of inequality (2.4), or inequality (2.4) proposed in this paper extends Nesbitt's inequality to general form.  

\subsection{Upper Bound of General Form of Nesbitt's  Inequality} 
Inequality (2.4) gives the lower bound. A very direct question is what is the upper bound? Obviously, generally speaking, the upper bound does not exist, or the upper bound can tend to infinity. This is easy to understand. For example, you can choose that ${{x}_{1}}$ is large and the rest of ${{x}_{i}}$ ($i\ne1$) is small, so that the left side of inequality (2.4) becomes very large and can be larger than any given large number. 

Of course, it is possible to obtain the upper bound of inequality (2.4) when ${{x}_{i}}$ is constrained by some other conditions besides being a positive number. As for the upper bound, I give the following theorem.

\begin{Theorem}\label{thm2.3}
	${{x}_{i}}\in{{\mathbf{R}}^{+}}$, $i=1,2,\cdots ,n$, and
	\begin{equation}\label{eq2.5}
		\sum\limits_{j=1}^{n}{{{x}_{j}}}-{{x}_{i}}>{{x}_{i}}
	\end{equation} then
	\begin{equation}\label{eq2.6}
		\sum\limits_{i=1}^{n}{\frac{{{x}_{i}}}{\sum\limits_{j=1}^{n}{{{x}_{j}}}-{{x}_{i}}}}<2
	\end{equation} 
\end{Theorem}

The proof of this theorem will be given later. Combining theorem 2.3 and theorem 2.2, the following theorem is obtained:

\begin{Theorem}\label{thm2.4}
	Let ${{x}_{i}}\in{{\mathbf{R}}^{+}}$, $i=1,2,\cdots ,n$, and $\sum\limits_{j=1}^{n}{{{x}_{j}}}-{{x}_{i}}>{{x}_{i}}$, then 
	\begin{equation}\label{eq2.7}
		\frac{n}{n-1}\le \sum\limits_{i=1}^{n}{\frac{{{x}_{i}}}{\sum\limits_{j=1}^{n}{{{x}_{j}}}-{{x}_{i}}}}<2
	\end{equation} 
	Equality holds iff ${{x}_{1}}={{x}_{2}}=\cdots ={{x}_{n}}$. 
\end{Theorem}

Nesbitt's inequality has been generalized by using  theorem 2.1 proposed above. Moreover, there are many other applications of theorem 2.1. Here are just a few examples. 

\section{Examples of Applications}
\paragraph{Example 3.1}\label{Example 3.1}
Let $n=3$, $S={{x}_{1}}+{{x}_{2}}+{{x}_{3}}$, by substituting the inequality  (2.1) of theorem 2.1, we obtain:
\[\begin{aligned}
 \frac{1}{{{x}_{2}}+{{x}_{3}}}+\frac{1}{{{x}_{3}}+{{x}_{1}}}+\frac{1}{{{x}_{1}}+{{x}_{2}}}&=\sum\limits_{i=1}^{3}{\frac{1}{S-{{x}_{i}}}}\ge \frac{3}{S-\bar{x}} \\ 
& =\frac{3}{{{x}_{1}}+{{x}_{2}}+{{x}_{3}}-\frac{1}{3}\left( {{x}_{1}}+{{x}_{2}}+{{x}_{3}} \right)} \\ 
& =\frac{9}{2\left( {{x}_{1}}+{{x}_{2}}+{{x}_{3}} \right)}  
\end{aligned}\]
Equality holds iff ${{x}_{1}}={{x}_{2}}={{x}_{3}}$. 

\paragraph{Example 3.2}\label{Example 3.2}
Let $n=3$, $S={{x}_{1}}+{{x}_{2}}+{{x}_{3}}$, ${{a}_{1}}={{x}_{2}}+{{x}_{3}}$, ${{a}_{2}}={{x}_{3}}+{{x}_{1}}$, ${{a}_{3}}={{x}_{1}}+{{x}_{2}}$, by substituting  inequality (2.1) of theorem 2.1, We get a well-known result as follows:

\[\begin{aligned}
\frac{1}{{{x}_{1}}}+\frac{1}{{{x}_{2}}}+\frac{1}{{{x}_{3}}}&=\sum\limits_{i=1}^{3}{\frac{1}{S-{{a}_{i}}}}\ge \frac{3}{S-\bar{a}} \\ 
& =\frac{3}{{{x}_{1}}+{{x}_{2}}+{{x}_{3}}-\frac{1}{3}\left( {{a}_{1}}+{{a}_{2}}+{{a}_{3}} \right)} \\ 
& =\frac{3}{{{x}_{1}}+{{x}_{2}}+{{x}_{3}}-\frac{2}{3}\left( {{x}_{1}}+{{x}_{2}}+{{x}_{3}} \right)} \\ 
& =\frac{9}{{{x}_{1}}+{{x}_{2}}+{{x}_{3}}}  
\end{aligned}\]
Equality holds iff ${{x}_{1}}={{x}_{2}}={{x}_{3}}$. 

This is a well-known result, but it can be obtained directly, conveniently and quickly by using theorem 2.1 of this paper. 

Because the conditions of the theorem (2.1) are easily satisfied, many inequalities can be easily constructed. Here are some examples, in the following examples, we assume that ${{x}_{i}}\in{{\mathbf{R}}^{+}}$. 
 
\paragraph{Example 3.3}\label{Example 3.3}
Let $n=3$, $S={{x}_{1}}+{{x}_{2}}+{{x}_{3}}$, by substituting inequality (2.1) of theorem 2.1, We have:
\[\begin{aligned}
\frac{1}{{{x}_{2}}+{{x}_{3}}}+\frac{1}{{{x}_{3}}+{{x}_{1}}}&=\sum\limits_{i=1}^{2}{\frac{1}{S-{{x}_{i}}}}\ge \frac{2}{S-\bar{x}} \\ 
& =\frac{2}{\left( {{x}_{1}}+{{x}_{2}}+{{x}_{3}} \right)-\frac{1}{2}\left( {{x}_{1}}+{{x}_{2}} \right)} \\ 
& =\frac{4}{2\left( {{x}_{1}}+{{x}_{2}}+{{x}_{3}} \right)-\left( {{x}_{1}}+{{x}_{2}} \right)} \\ 
& =\frac{4}{{{x}_{1}}+{{x}_{2}}+2{{x}_{3}}}  
\end{aligned}\] 
Equality holds iff ${{x}_{1}}={{x}_{2}}={{x}_{3}}$. 

\paragraph{Example 3.4}\label{Example 3.4}
Let $n=4$, $S={{x}_{1}}+{{x}_{2}}+{{x}_{3}}+{{x}_{4}}$, ${{a}_{1}}={{x}_{3}}+{{x}_{4}}$, ${{a}_{2}}={{x}_{2}}+{{x}_{4}}$, ${{a}_{3}}={{x}_{2}}+{{x}_{3}}$, ${{a}_{4}}={{x}_{1}}+{{x}_{4}}$, ${{a}_{5}}={{x}_{1}}+{{x}_{3}}$, ${{a}_{6}}={{x}_{1}}+{{x}_{2}}$, by substituting  inequality (2.1) of theorem 2.1, We have:
\[\begin{aligned}
& \frac{1}{{{x}_{1}}+{{x}_{2}}}+\frac{1}{{{x}_{1}}+{{x}_{3}}}+\frac{1}{{{x}_{1}}+{{x}_{4}}}+\frac{1}{{{x}_{2}}+{{x}_{3}}}+\frac{1}{{{x}_{2}}+{{x}_{4}}}+\frac{1}{{{x}_{3}}+{{x}_{4}}} \\ 
& =\sum\limits_{i=1}^{6}{\frac{1}{S-{{a}_{i}}}}\ge \frac{6}{S-\bar{a}} \\ 
& =\frac{6}{{{x}_{1}}+{{x}_{2}}+{{x}_{3}}+{{x}_{4}}-\frac{1}{6}\left( {{a}_{1}}+{{a}_{2}}+{{a}_{3}}+{{a}_{4}}+{{a}_{5}}+{{a}_{6}} \right)} \\ 
& =\frac{6}{{{x}_{1}}+{{x}_{2}}+{{x}_{3}}+{{x}_{4}}-\frac{3}{6}\left( {{x}_{1}}+{{x}_{2}}+{{x}_{3}}+{{x}_{4}} \right)} \\ 
& =\frac{12}{{{x}_{1}}+{{x}_{2}}+{{x}_{3}}+{{x}_{4}}} \\ 
\end{aligned}\] 
Equality holds iff ${{a}_{1}}={{a}_{2}}={{a}_{3}}={{a}_{4}}={{a}_{5}}={{a}_{6}}$, then we can deduce that the equality holds iff ${{x}_{1}}={{x}_{2}}={{x}_{3}}={{x}_{4}}$.

\paragraph{Example 3.5}\label{Example 3.5}
Let $n=4$, $S={{x}_{1}}+{{x}_{2}}+{{x}_{3}}+{{x}_{4}}$, by substituting inequality (2.1) of theorem 2.1, We have:
\[\begin{aligned}
& \frac{1}{{{x}_{2}}+{{x}_{3}}+{{x}_{4}}}+\frac{1}{{{x}_{3}}+{{x}_{4}}+{{x}_{1}}}+\frac{1}{{{x}_{4}}+{{x}_{1}}+{{x}_{2}}}+\frac{1}{{{x}_{1}}+{{x}_{2}}+{{x}_{3}}} \\ 
& =\sum\limits_{i=1}^{4}{\frac{1}{S-{{x}_{i}}}}\ge \frac{4}{S-\bar{x}} \\ 
& =\frac{4}{{{x}_{1}}+{{x}_{2}}+{{x}_{3}}+{{x}_{4}}-\frac{1}{4}\left( {{x}_{1}}+{{x}_{2}}+{{x}_{3}}+{{x}_{4}} \right)} \\ 
& =\frac{16}{3\left( {{x}_{1}}+{{x}_{2}}+{{x}_{3}}+{{x}_{4}} \right)} \\ 
\end{aligned}\] 
Equality holds iff ${{x}_{1}}={{x}_{2}}={{x}_{3}}={{x}_{4}}$.

\paragraph{Example 3.6}\label{Example 3.6}
Let $n=2$, $S={{\left( {{x}_{1}}+{{x}_{2}} \right)}^{2}}$, ${{a}_{1}}=x_{1}^{2}$, ${{a}_{2}}=x_{2}^{2}$, by substituting inequality (2.1) of theorem 2.1, We have:
\[\begin{aligned}
& \frac{1}{{{\left( {{x}_{1}}+{{x}_{2}} \right)}^{2}}-x_{1}^{2}}+\frac{1}{{{\left( {{x}_{1}}+{{x}_{2}} \right)}^{2}}-x_{2}^{2}}=\sum\limits_{i=1}^{2}{\frac{1}{S-{{a}_{i}}}}\ge \frac{2}{S-\bar{a}} \\ 
& =\frac{2}{{{\left( {{x}_{1}}+{{x}_{2}} \right)}^{2}}-\frac{1}{2}\left( {{a}_{1}}+{{a}_{2}} \right)}=\frac{2}{{{\left( {{x}_{1}}+{{x}_{2}} \right)}^{2}}-\frac{1}{2}\left( x_{1}^{2}+x_{2}^{2} \right)} \\ 
& =\frac{4}{{{\left( {{x}_{1}}+{{x}_{2}} \right)}^{2}}+2{{x}_{1}}{{x}_{2}}} \\ 
\end{aligned}\] 
i.e.
\[\frac{1}{\left( 2{{x}_{1}}+{{x}_{2}} \right){{x}_{2}}}+\frac{1}{\left( 2{{x}_{2}}+{{x}_{1}} \right){{x}_{1}}}\ge \frac{4}{{{\left( {{x}_{1}}+{{x}_{2}} \right)}^{2}}+2{{x}_{1}}{{x}_{2}}}\]
Equality holds iff ${{a}_{1}}={{a}_{2}}$, then we can deduce that the equality holds iff ${{x}_{1}}={{x}_{2}}$.

If let ${{x}_{1}}+{{x}_{2}}=1$, then the above inequality becomes
\[\frac{1}{1-x_{1}^{2}}+\frac{1}{1-x_{2}^{2}}\ge \frac{4}{1+2{{x}_{1}}{{x}_{2}}}\]
Equality holds iff ${{x}_{1}}={{x}_{2}}$.

\paragraph{Example 3.7 Application in Geometry}\label{Example 3.7}
The above results can be directly applied to geometry. Suppose ${{x}_{1}}$, ${{x}_{2}}$,…, ${{x}_{n}}$ are the lengths of the sides of a given polygon. Obviously, they satisfy the conditions of theorem 2.4. By using theorem 2.4, the following results can be obtained directly:
\[\begin{aligned}
& \frac{n}{n-1}\le \frac{{{x}_{1}}}{{{x}_{2}}+{{x}_{3}}+\cdots +{{x}_{n}}}+\frac{{{x}_{2}}}{{{x}_{3}}+{{x}_{4}}+\cdots +{{x}_{n}}+{{x}_{1}}} \\ 
& +\cdots +\frac{{{x}_{n}}}{{{x}_{1}}+{{x}_{2}}+\cdots +{{x}_{n-1}}}<2 \\ 
\end{aligned}\]
Equality holds iff ${{x}_{1}}={{x}_{2}}=\cdots ={{x}_{n}}$.

The above results can be directly applied to triangles. Assuming that ${{x}_{1}}$, ${{x}_{2}}$ and ${{x}_{3}}$ are the lengths of the sides of a given triangle, using theorem 2.4, let $n=3$, we can get the following corollary directly.

\begin{Corollary}\label{col3.1}
	Assuming that ${{x}_{1}}$, ${{x}_{2}}$ and ${{x}_{3}}$ are the lengths of the sides of a given triangle, we have:
	\[\frac{3}{2}\le \frac{{{x}_{1}}}{{{x}_{2}}+{{x}_{3}}}+\frac{{{x}_{2}}}{{{x}_{1}}+{{x}_{3}}}+\frac{{{x}_{3}}}{{{x}_{1}}+{{x}_{2}}}<2\]
	Equality holds iff ${{x}_{1}}={{x}_{2}}={{x}_{3}}$.
\end{Corollary}

This is a well-known result.
\section{Conclusions and Prospects}
The inequalities presented in this paper can solve a class of lower or upper bounds of fractional functions with $n$ variables such as $\sum\limits_{i=1}^{n}{\left( {1}/{S-{{x}_{i}}}\; \right)}$ or $\sum\limits_{i=1}^{n}{{{{x}_{i}}}/{\left( \sum\limits_{j=1}^{n}{{{x}_{j}}}-{{x}_{i}} \right)}\;}$. 

Because the condition $S-{{x}_{i}}>0$ ($i=1,2,\cdots ,n$) can be easily satisfied, many new inequalities can be easily created by using inequality (2,1) proposed in this paper. The generalization of Nesbitt's inequality to the general form of $n$ variables is a good example. 

The inequalities presented in this paper can be easily used to prove some known results, and the calculation is simple, fast, easy to remember and easy to use. 

Up to now, there are tens of thousands of inequalities (some of them are given in [9]). It is impossible for people to use them all, let alone memorize them, and there is no need to do so. 

The author believes that although it is impossible to deal with all inequalities in unique unified way, we should try our best to establish the internal relationship between seemingly different inequalities and adopt a unified way of thinking to deal with inequalities as much as possible, which is an important research field. The result of this paper is only the tip of the iceberg, and there is still much work to be done. 

I hope that the proposed theorems in this paper can be applied more widely.

\end{document}